\documentclass[10pt,reqno]{amsart}

\usepackage{amssymb}
\usepackage{amsthm}
\usepackage{amsmath}

\newtheorem{thm}{Theorem}

\newtheorem{cor}[thm]{Corollary}
\newtheorem{lem}[thm]{Lemma}

\newtheorem{knownthm}{Theorem}

\newtheorem{knownlem}[knownthm]{Lemma}
\theoremstyle{definition}
\newtheorem*{rem}{Remark}

\newcommand{\round}{\partial}

\newcommand{\C}{\mathbb{C}}
\newcommand{\D}{\mathbb{D}}

\newcommand{\R}{\mathbb{R}}

\newcommand{\A}{\mathcal{A}}
\renewcommand{\S}{\mathcal{S}}

\newcommand{\no}{\noindent}
\newcommand{\dstyle}{\displaystyle}
\renewcommand{\Re}{\textup{Re}\,}

\renewcommand{\qed}{\hfill $\square$}

\title[Boundedness for Robertson functions]{Boundedness, univalence and quasiconformal extension of Robertson functions}
\author[I. Hotta]{Ikkei Hotta}
\author[L.-M. Wang]{Li-Mei Wang}
\address{Division of Mathematics, Graduate School of Information Sciences, Tohoku University, 6-3-09 Aramaki-Aza-Aoba, Aoba-ku, Sendai, Miyagi 980-8579, Japan, phone : +81-22-795-4636}
\email{ikkeihotta@ims.is.tohoku.ac.jp}
\urladdr{http://www.ims.is.tohoku.ac.jp/{\textasciitilde}ikkeihotta}
\email{rime@ims.is.tohoku.ac.jp}
\date{\today}

\begin{document}

\maketitle

\begin{abstract}
This article contains several results for $\lambda$-Robertson functions, i.e., analytic functions $f$ defined on the unit disk $\D$ satisfying $f(0) = f'(0)-1=0$ and $\Re e^{-i\lambda}\{1+zf''(z)/f'(z)\} > 0$ in $\D$, where $\lambda \in (-\pi/2,\pi/2)$.
We will discuss about conditions for boundedness and quasiconformal extension of Robertson functions.
In the last section we provide another proof of univalence for Robertson functions by using the theory of L\"owner chains.
\end{abstract}

\section{Introduction}

\par Let $\mathcal{A}$ be the family of functions $f$ analytic in the unit disc $\mathbb{D}=\{z\in \mathbb{C}:\, |z|<1\}$ with the usual normalization $f(0)=f'(0)-1=0$, and $\S$ be the subclass of $\A$ consisting of functions univalent in $\D$.

Let $\lambda$ be a real constant between $-\pi/2$ and $\pi/2$. 
The curve $\gamma_{\lambda}(t)=\exp(te^{i\lambda})$, $t\in \mathbb{R}$, and its rotations $e^{i\theta}\gamma_{\lambda}(t)$, $\theta\in \mathbb{R}$, are called \textit{$\lambda$-spirals}. 
A domain $\Omega$ with $0\in \Omega$ is called \textit{$\lambda$-spirallike} (with respect to 0) if for every $w \in \Omega$, the $\lambda$-spiral which connects $w$ and 0 lies in $\Omega$. 
A function $f\in \mathcal{A}$ is said to be a \textit{$\lambda$-spirallike function} if $f$ maps $\mathbb{D}$ univalently onto a $\lambda$-spirallike domain and the class of such functions is denoted by $\mathcal{SP}(\lambda)$. 
Spirallike functions are introduced  by \v Spa\v cek \cite{Spacek:1932} in 1933. 
We note that $0$-spirallike functions are precisely starlike functions.

It is known that a necessary and sufficient condition for $f\in \mathcal{A}$ to be in $\mathcal{SP}(\lambda)$ is that
\[
\Re
\left\{
e^{-i\lambda}\frac{zf'(z)}{f(z)}
\right\}
>0
\]  
for all $z \in \D$.
In \cite{KimSugawa:pre01}, Kim and Sugawa introduce the notion of \textit{$\lambda$-argument}. 
Let us set $\theta=\arg_{\lambda} w$ if $w\in e^{i\theta}\gamma_{\lambda}(\mathbb{R})$. 
We note that $\arg_{0}w=\arg w$. For some more properties of $\lambda$-argument, the reader may be referred to \cite{KimSugawa:pre01}.
By utilizing $\lambda$-argument, another equivalence can be obtained
\[
f\in \mathcal{SP}(\lambda) \Leftrightarrow \frac{\partial}{\partial\theta}\left(\arg_{\lambda}f(re^{i\theta})\right)>0 \quad (\theta \in \R,\,0<r<1).
\]
For general references about spirallike functions, see e.g. \cite{Duren:1983} or \cite{AhujaSilverman:1991}.

A function $f\in \mathcal{A}$ is said to be a \textit{$\lambda$-Robertson function} \cite{Kulshrestha:1976} if $f$ satisfies
\[
\Re
\left\{
e^{-i\lambda}\left(1+\frac{zf''(z)}{f'(z)}\right)
\right\}
>0
\]
for all $z \in \D$. 
Let $\mathcal{R}(\lambda)$ be the set of those functions. 
The definition of $\lambda$-Robertson functions shows immediately that $\mathcal{R}(0)$ is precisely the class of convex functions which is usually denoted by $\mathcal{K}$. 
Furthermore in view of the definitions of spirallike and Robertson functions, for a function $f\in \mathcal{A}$ the following relations are true;
\begin{eqnarray}
\label{relationship}%%%%%%
f\in \mathcal{R}(\lambda) 
&\Leftrightarrow& 
zf'(z)\in \mathcal{SP}(\lambda)\\[5pt]
&\Leftrightarrow& 
\int_{0}^{z}f'(\zeta)^{\alpha}d\zeta\in \mathcal{K}\nonumber\\[5pt]
&\Leftrightarrow& 
\frac{\partial}{\partial\theta}\left[\arg_{\lambda}\left(\frac{\partial}{\partial\theta}f(re^{i\theta})\right)\right]>0 \quad (\theta\in\R,\,0<r<1),\nonumber
\end{eqnarray}
where $\alpha=e^{-i\lambda}/\cos\lambda$.
A distinguished member of $\mathcal{R}(\lambda)$ is 
\begin{equation}\label{extremal}
f_{\lambda}(z)=\frac{(1-z)^{1-2e^{i\lambda}\cos\lambda}-1}{2e^{i\lambda}\cos\lambda -1}.
\end{equation}

The class $\mathcal{R}(\lambda)$ was first introduced by Robertson \cite{Robertson:1969}.
He showed that all functions in $\mathcal{R}(\lambda)$ are univalent if $0 < \cos \lambda \leq x_{0}$, where $x_{0} =0.2034\cdots$ is the unique positive root of the equation $16x^{3} + 16 x^{2} + x -1 = 0$ (in the original paper $x_{0}$ is evaluated as $0.2315\cdots$ which seems to be erroneous \cite{KimSugawa:2007}).
Later Libera and Ziegler \cite{LiberaZiegler:1972} and Chichra \cite{Chichra:1975} gave some improvements on the range of $\lambda$ for which $\mathcal{R}(\lambda) \subset \S$ by estimating the norm of the Schwarzian derivatives for the class $\mathcal{R}(\lambda)$.
Finally Pfaltzgraff \cite{Pfaltzgraff:1975} showed that $\mathcal{R}(\lambda) \subset \S$ if $0 < \cos \lambda \leq 1/2$ or $\cos \lambda = 1$.
This value is best possible.
Indeed, Robertson also presented in \cite{Robertson:1969} a non-univalent function which belongs to $\mathcal{R}(\lambda)$ for each $\lambda$ in the range $1/2 < \cos \lambda < 1$ by making use of Roysters's example \cite{Royster:1965}
$f_{\mu}^{*}(z) = ((1-z)^{-\mu} - 1)/\mu$, where $\mu$ is a number which satisfies $\mu + 1 = |\mu + 1| e^{i\lambda}, |\mu| \leq 1, |\mu+1|>1$ and $|\mu-1|>1$.

The class of $\lambda$-Robertson functions has been investigated by various authors. 
Recently the class $\mathcal{R}(\lambda)$ is still an interesting topic in geometric function theory (e.g. \cite{PonnuYanagihara:2008}).
Actually,  under the relationship \eqref{relationship} many properties of Robertson functions follows from those of spirallike functions.
For instance the coefficient estimate of $\mathcal{R}(\lambda)$ is an easy consequence of a result of Zamorski \cite{Zamorski:1960} (see also \cite{BajpaiM:1973}).
For some more information about Robertson functions, the reader is referred to e.g. \cite[Section 8]{AhujaSilverman:1991}.

In the present paper we would like to give several new results for the $\lambda$-Robertson functions.
In section 2 we will show that $\lambda$-Robertson functions are bounded whenever $\cos\lambda<1/\sqrt{2}$ which improves a result of Kim and Sugawa in \cite{KimSugawa:2007}.
Quasiconformal extension criteria which are related with Robertson functions are shown in section 3.
One of the criteria is also obtained by using the technique of L\"owner's theory.
We will discuss this problem in the last section and give an explicit L\"owner chain for Robertson functions.

%BoundednessBoundednessBoundednessBoundednessBoundednessBoundedness
%BoundednessBoundednessBoundednessBoundednessBoundednessBoundedness
%BoundednessBoundednessBoundednessBoundednessBoundednessBoundedness
%BoundednessBoundednessBoundednessBoundednessBoundednessBoundedness

\section{Boundedness of $\mathcal{R}(\lambda)$}

%BoundednessBoundednessBoundednessBoundednessBoundednessBoundedness
%BoundednessBoundednessBoundednessBoundednessBoundednessBoundedness
%BoundednessBoundednessBoundednessBoundednessBoundednessBoundedness
%BoundednessBoundednessBoundednessBoundednessBoundednessBoundedness

%11111111111111111111111111111111111111111111111111111111111
%11111111111111111111111111111111111111111111111111111111111

\subsection{Result and auxiliary lemma}

%11111111111111111111111111111111111111111111111111111111111
%11111111111111111111111111111111111111111111111111111111111

The boundedness of $\lambda$-Robertson function is analyzed by Kim and Sugawa \cite{KimSugawa:2007}.
It can be stated as follows after being translated to our notations.

\begin{knownthm}[\cite{KimSugawa:2007}] 
$\lambda$-Robertson functions are bounded by a constant depending only on $\lambda$ when $\cos\lambda<1/2$.
\end{knownthm}

They remark that the bound $1/2$ cannot be replaced by any number greater than $1/\sqrt{2}$  since the function given by ($\ref{extremal}$) is unbounded when $\cos\lambda>1/\sqrt{2}$. Our next result will verify that the bound $1/\sqrt{2}$ is best possible.

\begin{thm}
Let $f\in \mathcal{R}(\lambda)$ with $\cos\lambda<1/\sqrt{2}$, then $f$ is bounded.
\end{thm}

In order to prove the above result, the growth theorem of spirallike functions in \cite{Singh:1969} or \cite{AhujaSilverman:1991} is needed. Since those known forms are complicated there, we simplify them as follows.

\begin{lem}\label{lem2}%%%%%%
Let $f\in \mathcal{SP}(\lambda)$, then for $|z|=r<1$, we have
\[
\Psi_{1}(r)\leq |f(z)|\leq \Psi_{2}(r)
\]
where
\[
\Psi_{1}(r)=\left|P_{\lambda}(re^{i\theta_{1}}) \right|=\frac{r\exp\left(-\sin2\lambda\arcsin(r\sin\lambda)\right)}{(r\cos\lambda-\sqrt{1-r^2\sin^2\lambda})^{2\cos^2\lambda}}
\]
and
\[
\Psi_{2}(r)=\left|P_{\lambda}(re^{i\theta_{2}}) \right|=\frac{r\exp\left(\sin2\lambda\arcsin(r\sin\lambda)\right)}{(r\cos\lambda-\sqrt{1-r^2\sin^2\lambda})^{2\cos^2\lambda}}
\]
where
\[P_{\lambda}(z)=\frac{z}{(1-z)^{1+e^{2i\lambda}}}
\]
belongs to $\mathcal{SP}(\lambda)$ and $\theta_{j}$ ($j=1,2$) satisfy
\[
\sin(\lambda+\theta_{j})=r\sin\lambda \hspace{15pt}(j=1,2)
\]
and $\cos(\lambda+\theta_{1})<0$ and $\cos(\lambda+\theta_{2})>0$ respectively. 
\end{lem}

%22222222222222222222222222222222222222222222222222222222222
%22222222222222222222222222222222222222222222222222222222222

\subsection{Proof of Theorem 1}

%22222222222222222222222222222222222222222222222222222222222
%22222222222222222222222222222222222222222222222222222222222

Equivalence (1) and Lemma \ref{lem2} show that 
\begin{eqnarray*}
|f(z)| &=&\left|\int_{0}^{z}f'(\zeta)d\zeta\right|=\left| \int_{0}^{r} \frac{z}{r}f'(tz/r)dt\right|\\[5pt]
&\leq &\int_{0}^{r}|f'(tz/r)|dt\leq \int_{0}^{r}\frac{\exp(\sin(2\lambda)\arcsin(t\sin\lambda))}{(\sqrt{1-t^2\sin^2\lambda}-t\cos\lambda)^{2\cos^2\lambda}}dt
\end{eqnarray*}
where $0<|z|=r<1$.

Since the numerator in the above integrand is bounded over $[0,1]$, it is sufficient to estimate only the denominator.

Upon a change in the variable $s=1-t$, we obtain
\begin{eqnarray*}
\sqrt{1-t^2\sin^2\lambda}-t\cos\lambda
&=&
\sqrt{1-(1-s)^2\sin^2\lambda}-(1-s)\cos\lambda\\
&=&
\sqrt{\cos^2\lambda+2s\sin^2\lambda-s^2\sin^2\lambda}-(1-s)\cos\lambda\\
&=&
\cos\lambda \sqrt{1+2s\tan^2\lambda-s^2\tan^2\lambda}-(1-s)\cos\lambda\\
&=&
\cos\lambda[1+1/2(2s\tan^2\lambda-s^2\tan^2\lambda)+O(s^2)]-(1-s)\cos\lambda\\
&=&
\frac{s}{\cos\lambda}+O(s^2)
\end{eqnarray*}
when $s\to 0$.

Therefore $f(z)$ is bounded whenever $2\cos^2\lambda<1$, that is, $\cos\lambda<1/\sqrt{2}$. The example given by $(\ref{extremal})$ ensures the sharpness of the value $1/\sqrt{2}$.
\qed

\begin{rem}
Note that our method is not applicable for the case when $\cos\lambda=1/\sqrt{2}$. Since the function $f_{\lambda}(z)$ given in ($\ref{extremal}$) is bounded when $\cos\lambda=1/\sqrt{2}$, we may expect that $\mathcal{R}(\lambda)$ consists of bounded functions as well in this case.
\end{rem}

%QCextensionQCextensionQCextensionQCextensionQCextensionQCextension
%QCextensionQCextensionQCextensionQCextensionQCextensionQCextension
%QCextensionQCextensionQCextensionQCextensionQCextensionQCextension
%QCextensionQCextensionQCextensionQCextensionQCextensionQCextension

\section{Quasiconformal Extension}

%QCextensionQCextensionQCextensionQCextensionQCextensionQCextension
%QCextensionQCextensionQCextensionQCextensionQCextensionQCextension
%QCextensionQCextensionQCextensionQCextensionQCextensionQCextension
%QCextensionQCextensionQCextensionQCextensionQCextensionQCextension

%11111111111111111111111111111111111111111111111111111111111
%11111111111111111111111111111111111111111111111111111111111

\subsection{Results}

%11111111111111111111111111111111111111111111111111111111111
%11111111111111111111111111111111111111111111111111111111111

In this section we would like to discuss about the new quasiconformal extension criteria for Robertson functions.
Let us denote by $\S(k)$ the family of functions lie in $\S$ and can be extended to quasiconformal automorphisms of $\C$ so that the complex dilatation $\mu_{f} = (\round f / \round \bar{z}) / (\round f / \round z)$ satisfies $|\mu_{f}(z)| \leq k < 1$ for almost every $z \in \C$.

We will show the following which is an extension of a result of Ruscheweyh \cite[Corollary 1]{Ruscheweyh:1976};

\begin{thm}\label{result}%%%%%%
Let $f \in \A,\,k \in [0,1)$ and $\lambda\in(-\pi/2,\pi/2),\,q>-1$ be related by 
\begin{equation}\label{eq02}%%%%%%%
0 < \cos \lambda \leq \left\{
\begin{array}{llc}
k/2, & \textit{if}& -1 < q \leq 0 , \\[5pt]
k/(2+4q), & \textit{if}&  0 < q.
\end{array}
\right.
\end{equation}
If f satisfies
\begin{equation}\label{eq01}%%%%%%%%
\Re
\left\{
e^{-i\lambda}
\left(
1 + \frac{zf''(z)}{f'(z)} + q \frac{zf'(z)}{f(z)}
\right)
\right\}
>0
\end{equation}
for all $z \in \D$, then $f \in \S(k)$.
If, in addition, $f''(0)=0$, $\eqref{eq02}$ can be replaced by 
\begin{equation}\label{eq05}%%%%%%%
0 < \cos \lambda \leq \left\{
\begin{array}{llc}
k, & \textit{if}& -1 < q \leq 0 , \\[5pt]
k/(1+2q), & \textit{if}&  0 < q.
\end{array}
\right.
\end{equation}

\end{thm}

We note that when $q=0$ Theorem 3 claims quasiconformal extension of $\lambda$-Robertson functions which can be stated explicitly as follows;

\begin{cor}\label{qccor}
Let $f \in \mathcal{R}(\lambda)$ with $\lambda\in (-\pi/2,\pi/2)$ satisfying
\begin{equation*}\label{eq06}%%%%%%
0 < \cos \lambda \leq k/2,
\end{equation*}
then $f\in\S(k)$. 
If, in addition, $f''(0)=0$ and $\eqref{eq06}$ can be replced by
\begin{equation*}
0 < \cos \lambda \leq k,
\end{equation*}
 then $f\in \S(k)$.
\end{cor}

\par We note here that the second case in Corollary 4 also implies that function $f \in \mathcal{R}(\lambda)$ with $f''(0)=0$ for arbitrary $\lambda\in (-\pi/2,\pi/2)$ is univalent which was proved by Singh and Chichra \cite{SinghChichra:1977b} by means of Ahlfors's criterion for univalence as well.

\subsection{Preliminaries}
The following several results will be used later in our arguments.
Here, set
$$
H_{s}(z) = s \left(1 + \frac{zf''(z)}{f'(z)}\right) + (1-s)\frac{zf'(z)}{f(z)}.
$$

\begin{knownthm}[\cite{Hotta:2010b}]\label{main}%%
Let $a>0,\,b \in \R,\,s = a + i b,\,k \in [0,1)$ and $f \in \A$.
Assume that for a constant $c \in \C$ and all $z \in \D$
\begin{equation*}
\left| 
c|z|^{2}+s-a(1-|z|^{2})H_{s}(z)
\right|
\leq M
\end{equation*}
with
$$
M =
\left\{
\begin{array}{ll}
a k |s| + (a-1)|s+c|, &\textit{if}\hspace{10pt}  0 < a \leq 1, \\[5pt]
k |s|, &\textit{if}\hspace{10pt} a>1,
\end{array} 
\right.
$$
then $f \in \S(l)$, where
\begin{equation*}\label{mainl}%%
l=
\frac
{
2ka+(1-k^{2}) |b|
}
{
(1+k^{2})a+(1-k^{2}) |s|
}
<1.
\end{equation*}
\end{knownthm}

We note that in the above theorem $l=k$ if and only if $b=0$ (\cite{Hotta:2010b}).

\begin{knownlem}[e.g. \cite{Ruscheweyh:1976}]\label{lemB}%%%%%%%
Let $p(z) = 1 + a_{n}z^{n} + \cdots$ be analytic and $\Re p(z)>0$ on $\D$.
Then
\begin{equation*}
\left|
p(z) -1-\frac{2|z|^{2n}}{1-|z|^{2n}}
\right|
\leq 
\frac{2|z|^{2n}}{1-|z|^{2n}}
\end{equation*}
for all $z \in \D$.
\end{knownlem}

%22222222222222222222222222222222222222222222222222222
%22222222222222222222222222222222222222222222222222222

\subsection{Proof of Theorem \ref{result}}

%22222222222222222222222222222222222222222222222222222
%22222222222222222222222222222222222222222222222222222

Let $s = 1/(1+q)$ and $f(z)=z+\sum_{n=2}^{\infty}a_{n}z^n$, then for 
\begin{eqnarray*}
p(z) &=& \frac{e^{-i\lambda} H_{s}(z) + i\sin \lambda}{\cos \lambda}\\
&=& 1 + \frac{e^{-i\lambda}}{\cos \lambda} (s+1)a_{2}z + \cdots.
\end{eqnarray*}
we have $p'(0)=0$ if and only if $f''(0)=0$.
Condition \eqref{eq01} implies that $p(z)$ is analytic in $\D$ and fulfills $\Re p(z) > 0$ for all $z \in \D$.
With $\dstyle (c+s) = \frac 2 n se^{i\lambda}\cos\gamma,\,n=1, 2$, we obtain from Lemma $\ref{lemB}$ that
$$
\begin{tabular}{llllll}
\multicolumn{1}{l}
{$\dstyle \left|\frac{(c+s)|z|^{2}}{1-|z|^{2}} - s(H_{s}(z) -1)\right|$} \\[12pt]
$\hspace{60pt}\dstyle \leq
s|\cos \lambda|
\left\{
\left|\frac{2|z|^{2n}}{1-|z|^{2n}} - (p(z) -1)\right| + 
\left|\frac{2|z|^{2n}}{1-|z|^{2n}} - \frac2n \frac{|z|^{2}}{1-|z|^{2}}\right|
\right\}$\\[12pt]
$\hspace{60pt}\dstyle \leq \frac{2s}{n}\frac{|\cos \lambda|}{1-|z|^{2}}$.
\end{tabular}
$$

\no
Therefore by Theorem \ref{main} $f$ can be extended to a $k$-quasiconformal automorphism of $\C$ whenever
$$
\frac2n s|\cos \lambda| \leq \left\{
\begin{tabular}{llc}
$\dstyle k s^{2} + \frac2n s |\cos \lambda| (s-1)$, & \textit{if}&$0 < s \leq 1$, \\[5pt]
$k s$, & \textit{if}& $1 < s$,
\end{tabular}
\right.
$$
which is equivalent to \eqref{eq02} if $n=1$ and to \eqref{eq05} if $n=2$.
\qed

%LoewnerchainsLoewnerchainsLoewnerchainsLoewnerchainsLoewnerchainsLoewnerchains
%LoewnerchainsLoewnerchainsLoewnerchainsLoewnerchainsLoewnerchainsLoewnerchains
%LoewnerchainsLoewnerchainsLoewnerchainsLoewnerchainsLoewnerchainsLoewnerchains
%LoewnerchainsLoewnerchainsLoewnerchainsLoewnerchainsLoewnerchainsLoewnerchains

\section{L\"owner chain}

%LoewnerchainsLoewnerchainsLoewnerchainsLoewnerchainsLoewnerchainsLoewnerchains
%LoewnerchainsLoewnerchainsLoewnerchainsLoewnerchainsLoewnerchainsLoewnerchains
%LoewnerchainsLoewnerchainsLoewnerchainsLoewnerchainsLoewnerchainsLoewnerchains
%LoewnerchainsLoewnerchainsLoewnerchainsLoewnerchainsLoewnerchainsLoewnerchains

We can find another proof for univalency of Robertson functions by making use of the theory of L\"owner chains. 

The following theorem is well known.
Here, we denote $\round f / \round t$ and $\round f / \round z$ by $\dot{f}$ and $f'$ respectively.

\begin{knownthm}[\cite{Pom:1965}, see also \cite{Hotta:2010a}]
Let $0 < r_{0} \leq 1$.
Let $f_{t}(z) = \sum_{n=1}^{\infty}a_{n}(t)z^{n}$, $a_{1}(t) \neq 0$,\, be analytic for each $t \in [0,\infty)$ in $\D_{r_{0}}$ and locally absolutely continuous in $[0,\infty)$, locally uniformly with respect to $\D_{r_{0}}$, where $a_{1}(t)$ is a complex-valued function on $[0,\infty)$.
For almost all $t \in [0,\infty)$ suppose
\begin{equation}\label{LDE}%%%%%%
\dot{f_{t}}(z) =z f_{t}'(z) p_{t}(z)  \hspace{20pt} (z \in \D_{r_{0}})
\end{equation}
where $p_{t}(z)$ is analytic in $\D$ and satisfies $\Re p_{t}(z)>0,\,z \in \D$.
If $|a_{1}(t)| \to \infty$ for $t \to \infty$ and if $\{f_{t}(z)/a_{1}(t)\}$ forms a normal family in $\D_{r_{0}}$, then for each $t \in [0,\infty)$ $f_{t}(z)$ can be continued analytically to a univalent function on $\D$.
\end{knownthm}

The next lemma is needed for our discussion;

\begin{knownlem}[\cite{Wang:preprint}, Theorem 3]\label{lemma3}%%%%%%
Suppose that $\lambda \in (-\pi/2,\pi/2)$.
Let $p(z)$ be an analytic function defined on $\D$ which satisfies $p(0)=1$ and $\Re e^{-i\lambda}p(z) >0$ for all $z \in \D$.
Then we have
\begin{equation*}
\left|p(z)-\left(\frac{1}{1-r^2}+\frac{r^2}{1-r^2}e^{2i\lambda}\right)\right|\leq 
\frac{2r}{1-r^2}\cos \lambda.
\end{equation*}
where $r=|z|<1$.

\end{knownlem}

Now we suppose that $|\lambda| \in [\pi/3, \pi/2)$ and $f$ is a $\lambda$-Robertson function.
Let us put
\begin{equation}\label{LC}%%%%%%
f_{t}(z) = f(e^{-t}z) - e^{-2i\lambda}(e^{2t}-1)e^{-t}zf'(e^{-t}z).
\end{equation}
Here we should note that a more general form of \eqref{LC} appears in \cite{Ruscheweyh:1976}.

It suffices to prove that $p_{t}(z)$ in \eqref{LDE} lies in the right-hand side of the complex plane $\C$
for all $z \in \D$ and a.e. $t \in [0,\infty)$. 
This is equivalent to
$$
\left|
\frac{\dot{f_{t}}(z) - zf_{t}'(z)}{\dot{f_{t}}(z) +zf_{t}'(z)}
\right|
< 1.
$$
Then a calculation shows that
\begin{equation}\label{eq03}%%
\left|
e^{-2t} e^{2i\lambda} + 1 - 
\left(
1 - e^{-2t}
\right)
\left(
1 + \frac{e^{-t}zf''(e^{-t}z)}{f'(e^{-t}z)}
\right)
\right|
<1
\end{equation}
implies univalency of $f$ and \eqref{eq03} follows from maximum modulus principle and Lemma \ref{lemma3} when $\cos \lambda \leq 1/2$.

\begin{rem}
Applying Becker's theorem \cite{Becker:1976} with \eqref{LC}, we also obtain the quasiconformal extension criterion for $\mathcal{R}(\lambda)$ which is in Corollary \ref{qccor}.
\end{rem}

\section*{Acknowledgement}
The authors would like to express their deep gratitude to Professor Toshiyuki Sugawa. 
This work would not finish without his useful discussion and constant encouragement.

%ReferencesReferencesReferencesReferencesReferencesReferencesReferences
%ReferencesReferencesReferencesReferencesReferencesReferencesReferences
%ReferencesReferencesReferencesReferencesReferencesReferencesReferences
%ReferencesReferencesReferencesReferencesReferencesReferencesReferences

\bibliographystyle{amsplain}

\end{document}